\documentclass[11pt,leqno]{article}
\usepackage{amsthm,amsfonts,amssymb,epsfig,graphics,amsmath,eufrak}
 \usepackage[latin1]{inputenc}\relax

\newcommand{\CalG}{\mathcal{G}}
\newcommand{\CalP}{\mathcal{P}}

\newcommand{\RR}{{\mathbb R}}

\newcommand{\EE }{{\mathbb E}}

\newcommand\cP{{\cal  P}}

\newcommand{\bC}{\mathbb{C}}

\newcommand\adots{\mathinner{\mkern2mu\raise1pt\hbox{.}
\mkern3mu\raise4pt\hbox{.}\mkern1mu\raise7pt\hbox{.}}}

\newcommand{\rank}{{\rm rank }}

\newtheorem{theo}{Theorem}[section]
\newtheorem{prop}[theo]{Proposition}

\newtheorem{lem}[theo]{Lemma}
\newtheorem{defi}[theo]{Definition}
\newtheorem{ass}[theo]{Assumption}

\newtheorem{exam}[theo]{Example}
\newtheorem{rem}[theo]{Remark}

\numberwithin{equation}{section}

\title{
Uniform stability estimates for constant-coefficient symmetric
hyperbolic boundary value problems }

\author{Olivier Gu\`es\footnote{Universit\'e de Provence, partially supported by European network HYKE,  HPRN-CT-2002-00282},
Guy M\'etivier\footnote{Universit\'e de Bordeaux, partially
supported by European network HYKE,  HPRN-CT-2002-00282.}, Mark
Williams\footnote{University of North Carolina, partially
supported  by NSF grants DMS-0070684 and DMS-0401252},
 Kevin Zumbrun\footnote{Indiana University, partially supported  by NSF grants  DMS-0070765 and DMS-0300487; K.Z. thanks Universit\'e de
Provence for its hospitality during the month-long visit in which
this work was carried out.
 }}

\date{August 21, 2005}

\begin{document}

\maketitle
\begin{abstract}
Answering a question left open in \cite{MZ2}, we show for general
symmetric hyperbolic boundary problems with constant coefficients,
including in particular systems with characteristics of variable
multiplicity, that the uniform Lopatinski condition implies strong
$L^2$ well-posedness, with no further structural assumptions. The
result applies, more generally, to any system that is strongly
$L^2$ well-posed for at least one boundary condition. The proof is
completely elementary, avoiding reference to Kreiss symmetrizers
or other specific techniques. On the other hand, it is specific to
the constant-coefficient case; at least, it does not translate in
an obvious way to the variable-coefficient case.  The result in
the hyperbolic case is derived from a more general principle that
can be applied, for example, to parabolic or partially parabolic
problems like the Navier-Stokes or viscous MHD equations
linearized about a constant state or even a viscous shock.
\end{abstract}

\bigbreak
\section{Introduction}\label{s1}
Consider a noncharacteristic, hyperbolic boundary value problem
with constant coefficients on the half-space
$\mathbb{R}^{d+1}_+=\{(t,x):x_d\geq 0\}$:
\begin{align}\label{hyp}
\begin{split}
&(a)\;Lu:=u_t + \sum_{j=1}^d A^j u_{x_j}= f\\
&(b)\;\Gamma u(t,\tilde x, 0)=g,
\end{split}
\end{align}
where $u\in \RR^n$, $\det A_d\ne 0$, $\tilde x:=(x_1, \dots,
x_{d-1})$, $\Gamma$ is a constant $k\times n$ matrix, and the
symbol $\sum_j^d A^ji \xi_j$ satisfies the hyperbolicity condition
\begin{equation}\label{hypcond}
 \sum_j^d A^ji \xi_j \, \text{\rm has only
 pure imaginary, semisimple eigenvalues for all } \xi\in \RR^d.
\end{equation}

There are two distinct, but partially overlapping classes of
systems for which the existence/stability theory is well
developed, namely the Friedrichs symmetrizable hyperbolic systems
with maximally dissipative boundary conditions  and the
Kreiss--M\'etivier class of strictly hyperbolic or
constant-multiplicity systems with $\Gamma$ satisfying a sharp
spectral condition called the uniform Lopatinski condition.

\begin{defi}\label{a1}
1.  The operator $L$ \eqref{hyp}(a) is called \emph{Friedrichs
symmetrizable} when there exists a positive symmetric matrix $S$
such that $SA^j$ is symmetric for $j=1,\dots,d$.

2.  Suppose $L$ is Friedrichs symmetrizable with symmetrizer $S$.
The boundary condition is maximally dissipative when $\rank
\;\Gamma=k$, $SA_d$ is negative definite on $\ker\Gamma$, and
$k=$dimension of the unstable subspace of $A_d$.

\end{defi}

\begin{rem}\label{aa0}
Let us recall a few well-known properties of the systems just
defined (see, e.g., \cite{Met4}, Chapter 2).

1.  $SA_d$ is negative definite on $\ker\Gamma$ if and only if
there are positive constants $c$ and $C$ such that for all
$h\in\bC^n$
\begin{align}\label{a01}
-(SA_dh,h)\geq c|h|^2-C|\Gamma h|^2.
\end{align}

2. One can define an adjoint problem $(L^*,\Gamma^*)$ where
\begin{align}\label{a02}
\begin{split}
&L^*=-\partial_t-\sum^d_1A_j^*\partial_{x_j},\\
&\Gamma^*\text{ is an }(n-k)\times n\text{ matrix with
}\ker\Gamma^*=(A_d\ker\Gamma)^\perp.
\end{split}
\end{align}
The problem $(L^*,\Gamma^*)$ is symmetrizable and maximally
dissipative in the backward sense; that is, $S^{-1}$ is a
symmetrizer for $-L^*$ and $S^{-1}A_d^*$ is positive definite on
$\ker \Gamma^*$.

3. Given a Friedrichs symmetrizable operator $L$ \eqref{hyp}(a)
one can always define a maximally dissipative boundary condition
for it using the projector $\pi_+$ of $\mathbb{C}^n$ onto the
unstable subspace $\mathcal{U}$ of $SA_d$. More precisely, if
$\dim \mathcal{U}=n_+$, one can take $\Gamma=T\pi_+$, where $T$ is
linear isomorphism
\begin{align}\label{at}
T:\mathcal{U}\to\mathbb{C}^{n_+}.
\end{align}
 We invite the
reader to check \eqref{a01} in that case.

\end{rem}

We will also consider more general boundary conditions of the form
\begin{align}\label{aa1}
\Gamma_\gamma u:=e^{\gamma t}\Gamma(D_t,D_{\tilde
x},\gamma)e^{-\gamma t}u=g
\end{align}
where $\Gamma(D_t,D_{\tilde x},\gamma)$ is a  Fourier multiplier:
\begin{align}\label{aa2}
\widehat{\Gamma
v}(\tau,\eta):=\Gamma(\tau,\eta,\gamma)\hat{v}(\tau,\eta)
\end{align}
defined by a continuous bounded $k\times n$ symbol
$\Gamma(\tau,\eta,\gamma)$.

As described in \cite{BT, MZ2}, physical applications such as
shock stability in magnetohydrodynamics (MHD) motivate the study
of a third class consisting of symmetric hyperbolic problems with
uniform Lopatinski boundary conditions but possibly
variable-multiplicity characteristics. This class was treated in
depth in \cite{MZ2} under some additional structural assumptions
on the system, satisfied in particular for MHD, at both the
linearized (constant- and variable-coefficient) and nonlinear
level, using a generalization of the symmetrizer techniques
introduced by Kreiss \cite{K} in the strictly hyperbolic setting.
However, it was noted that these structural assumptions could be
significantly relaxed in the constant-coefficient case for which
symmetrizers need not be smooth. Indeed, the construction in this
case hints of further generality, suggesting that for Friedrichs
symmetrizable systems, the uniform Lopatinski condition alone is
perhaps all that is needed for $L^2$ well-posedness (Definition
\ref{aa4}).

The purpose of this note is to verify by a very simple argument,
bypassing completely the symmetrizer constructions of \cite{K,
Met3, MZ2} that this conjecture is indeed correct. However, the
argument does not, at least in an obvious fashion, carry through
to the variable-coefficient or nonlinear case, for which the
Kreiss symmetrizer approach remains up to now the only choice.

\begin{defi}\label{aa4}  We say that the problem $(L,\Gamma_\gamma)$  \eqref{hyp}, \eqref{aa1} is
\emph{strongly $L^2$ well-posed} if there exists a $C>0$ such that
for $\gamma>0$, $f\in e^{\gamma t}L^2(\mathbb{R}^{d+1}_+)$, $g\in
e^{\gamma t}L^2(\mathbb{R}^d)$ there exists a unique solution
$u\in e^{\gamma t}L^2(\mathbb{R}^{d+1}_+)$, and $u$ satisfies the
energy estimate
\begin{equation}\label{uni}
\begin{aligned}
\gamma \int_{-\infty}^\infty &e^{-2\gamma t}\|u(\cdot,
t)\|_{L^2}^2 \, dt +
\int_{-\infty}^\infty e^{-2\gamma t} |u(0,t)|^2 \, dt \le\\
&\quad C\left(\gamma^{-1}\int_{-\infty}^\infty e^{-2\gamma
t}\|f(\cdot, t)\|_{L^2}^2 \, dt
 +
\int_{-\infty}^\infty  e^{-2\gamma t}|g(0,t)|^2 \, dt\right).
\end{aligned}
\end{equation}
The word \emph{strongly} is used to highlight the trace estimate
of $u$.
\end{defi}


 For Friedrichs symmetric systems with maximally dissipative
boundary conditions, strong $L^2$ well-posedness follows by
standard arguments (see, e.g., \cite{Met4}, Chapter 2) from an a
priori estimate of the form \eqref{uni} for the original problem
$(L,\Gamma)$ and an analogous estimate for the adjoint problem
$(L^*,\Gamma^*)$. The forward estimate, for example, is obtained
using integration by parts, taking the $L^2$ inner product of $Su$
against equation \eqref{hyp}; maximal dissipativity of $\Gamma$
\eqref{a01} allows the resulting boundary term to be estimated in
a straightforward way, yielding \eqref{uni}.   The adjoint
estimate is similar.

Maximally dissipative boundary conditions are clearly quite
special.  In order to define boundary conditions satisfying the
more general uniform Lopatinski condition, we first apply to
\eqref{hyp},\eqref{aa1} the Laplace transform in the temporal
variable $t$ and the Fourier transform in tangential spatial
variables $\tilde x:=(x_1, \dots, x_{d-1})$ to obtain the
resolvent equation:
\begin{equation}\label{ieq}
\begin{aligned}
\hat u'-G(\Lambda)\hat u&=\hat f,\\
\Gamma(\Lambda) \hat u(0)&= \hat g.
\end{aligned}
\end{equation}
Here $\hat u$, $\hat f$, and $\hat g$ denote the Laplace--Fourier
transforms of $u$, $f$, and $g$,
\begin{align}\label{aa3}
 \Lambda=(\tau,
\eta,\gamma)\in\cP:=\{(\tau,\eta,\gamma):(\tau,\eta)\in\mathbb{R}^d,\gamma>0\},
\end{align}
and
\begin{equation}\label{Geq}
G(\Lambda):= -A_d^{-1}\left(\gamma + i\tau +
i\sum_{j=1}^{d-1}\eta_j A^j\right).
\end{equation}
Recall that the Laplace transform of a function $f(t)\in e^{\gamma
t}L^2(t)$ is the Fourier transform of $e^{-\gamma t}f$.

From hyperbolicity, \eqref{hypcond}, we find easily the result of
Hersch \cite{H} that, for $\gamma>0$, $G(\Lambda)$ has no pure
imaginary eigenvalues. For, existence of an eigenvalue $i\kappa$,
$\kappa\in \RR$ of $G$ would imply existence of an eigenvalue
$\gamma + i\tau$ with nonzero real part $\gamma$ of the matrix
symbol $\sum_{j=1}^d A^j i\xi_j$, where $\xi=(\tilde \xi,
\xi_d):=(\eta, -\kappa) \in \RR^d$. Thus, defining
$\EE_-(\Lambda)$ to be the stable subspace of $G(\Lambda)$, we
have that $\dim \EE_-(\Lambda)$ is constant for all $\gamma>0$,
and (taking $\eta$, $\tau=0$, $\gamma=1$)
\begin{equation}\label{Edim}
\dim \EE_-(\Lambda) \equiv n_+ \quad \text{\rm for } \gamma>0,
\end{equation}
where $n_+$ denotes the
dimension of the unstable eigenspace of $A^d$.


\begin{defi}\label{a3}
A system $(L,\Gamma_\gamma)$ \eqref{hyp}, \eqref{hypcond},
\eqref{aa1} is said to satisfy the \emph{uniform Lopatinski
condition} when

\begin{align}\label{lopcond}
\begin{split}
&(i)k=\rank\;\Gamma(\Lambda)=\dim E_-(\Lambda)\text{ for all }\Lambda\in\mathcal{P}\\
&(ii) |v|\le C|\Gamma(\Lambda) v| \,\text{\rm for } v\in
\EE_-(\Lambda), \text{ for }C>0 \text{ independent of }
\Lambda\in\mathcal{P}.
\end{split}
\end{align}

\end{defi}

Kreiss \cite{K} showed that the uniform Lopatinski condition can
be derived as a necessary condition for strong $L^2$
well-posedness. The existence part of Definition \ref{aa4} applied
to the transformed problem \eqref{ieq} implies
\begin{align}\label{a3a}
\dim E_-(\Lambda)\geq k,
\end{align}
since, when $\hat{f}=0$, solutions of \eqref{ieq} in $L^2(x_d)$
must have boundary data in $E_-(\Lambda)$.  Plancherel's theorem
yields an estimate similar to \eqref{uni} for the transformed
problem (with, e.g., $L^2(t,\tilde{x})$ norms replaced by
$L^2(\tau,\eta)$ norms).  In fact, by studying special solutions
of \eqref{hyp}, \eqref{aa1} built from plane waves, this can be
pushed further to obtain estimates for \eqref{ieq} uniform with
respect to $\Lambda$:
\begin{align}\label{a3b}
\gamma\|\hat u\|^2_{L^2(x_d)}+|\hat u(0)|^2\leq C(\|\hat
f\|^2_{L^2(x_d)}/\gamma+|\hat g|^2)
\end{align}
for $\gamma>0$ and $C$ independent of $\Lambda\in\mathcal{P}$,
where $|\cdot|$ is the standard complex modulus (e.g., see
\cite{Met4}, Prop. 6.2.2).   Taking $\hat f=0$ in \eqref{a3b} we
deduce
\begin{align}\notag
|v|\leq C|\Gamma(\Lambda)|\text{ for }v\in E_-(\Lambda),
\end{align}
which implies
\begin{align}\label{a3c}
\dim E_-(\Lambda)\leq \rank\;\Gamma(\Lambda)\leq k\text{ for all
}\Lambda\in\mathcal{P}.
\end{align}
With \eqref{a3a} this shows the uniform Lopatinski condition is
necessary for strong $L^2$ well-posedness.

A major contribution of Kreiss \cite{K} was to show, in the
strictly hyperbolic case, by an ingenious construction of
frequency-dependent symmetrizers, that the uniform Lopatinski
condition is in fact equivalent to strong $L^2$ well-posedness, a
result later generalized to constant multiplicity hyperbolic
systems through the work of Majda--Osher \cite{MO} and Metivier
\cite{Met2}, and to certain variable-multiplicity hyperbolic
systems in \cite{MZ2}.

Our main result is the following extension to general Friedrichs
symmetrizable systems in the constant coefficient case:
\begin{theo}\label{hypmain}
Consider a constant coefficient Friedrichs symmetrizable  system
$L$ \eqref{hyp}(a) with boundary condition $\Gamma_\gamma$
\eqref{aa1}.   The system $(L,\Gamma_\gamma)$ is strongly $L^2$
well-posed if and only if it satisfies the uniform Lopatinski
condition.
\end{theo}

\begin{rem}\label{new}

1.  For constant coefficient systems one might try to obtain the
Kreiss estimate \eqref{uni} by direct estimation of solutions
constructed by Fourier-Laplace transform.  As far as we know this
has been done successfully only under more restrictive hypotheses
than the ones we make here (restrictions on multiplicities, order
of glancing points, etc.).
Weaker bounds (Hadamard well-posedness: in effect,
estimates exhibiting a loss of several derivatives)
have been established by this approach in great generality \cite{H}.

2.  For the constant coefficient systems we consider here (Friedrichs
symmetrizable with uniform Lopatinski boundary conditions), the
Kreiss estimate \eqref{uni} has been obtained by a simple
integration by parts argument when $f=0$ (\cite{S}, p. 199).
However, that argument does not appear to extend to the case $f\neq 0$.
On the other hand, an estimate losing one-half derivative
may easily be obtained by subtracting out the solution $w$
of the Cauchy problem extended to the whole space and
solving the residual problem with zero interior data and
boundary data $g-\Gamma w(0)$,
controlling $|w(0)|$ by the standard trace estimate
$|w(0)|\le C_1|w|_{H^{1/2}}\le C_2|f|_{H^{1/2}}$.

3.   Theorem \ref{main} was established using symmetrizers in
\cite{MZ2} under the additional structural assumption that, at
frequencies $\xi_0$ in the vicinity of which the eigenvalues
$a_j(\xi)$ of the symbol $\sum_{j=1}^d A^ji\xi_j$ are of variable
multiplicity, crossing eigenvalues are either {\it geometrically
regular} in the sense that eigenvalues and eigenprojections are
both analytic, {\it totally nonglancing} in the sense that
$\partial a_j/\partial \xi_d$ have a common, nonzero sign for all
$a_j$ involved, or {\it linearly separating} in the sense that
crossing eigenvalues $a_j(\xi)$ separate to linear order in the
distance of $\xi$ from a smooth manifold where they agree. The new
content of Theorem \ref{hypmain} is that these additional
assumptions may be dropped.

\end{rem}

As sketched briefly in Section \ref{viscous}, the same argument
yields an analogous result for the linearized equations about a
planar viscous shock or boundary layer with ``real'', or physical,
viscosity.
Thus, the general principle contained in Proposition \ref{main}
can be also be applied to parabolic or partially parabolic
problems.  However, as discussed in Section \ref{var}, our results
do not apply to the nonlinear or variable-coefficient case, either
in the hyperbolic or viscous--hyperbolic context.

\bigbreak
\section{Generalized resolvent-type equations}\label{general}


It remains to prove the sufficiency of the uniform Lopatinski
condition in Theorem \ref{hypmain}.  We'll deduce this from the
theory of maximally dissipative problems together with a new
result, Proposition \ref{main}, for constant-coefficient
``generalized resolvent-type'' equations
\begin{align}\label{eq}
\begin{split}
&L(\Lambda)u:=u'-G(\Lambda)u=f,\\
&\Gamma(\Lambda) u(0)=g.
\end{split}
\end{align}
on the half-line $x\in \RR^+$.  Here $\Lambda$ is a parameter
confined to a connected open set $\CalP$, and $\Gamma(\Lambda)$ is
a $k\times n$ matrix.  Initially, the \emph{only} assumption we
make about the $n\times n$ matrix $G(\Lambda)$ is that it has no
pure imaginary eigenvalues for $\Lambda\in\mathcal{P}$.  The
parameter $\Lambda$ might represent Laplace and/or Fourier
frequencies, model variables, etc.. If we define $\EE_-(\Lambda)$
to be the stable subspace of $G(\Lambda)$, these hypotheses imply
that $\dim \EE_-(\Lambda)$ is independent of $\Lambda\in\CalP$.

\begin{defi}\label{stab}
Relative to some choice of $\alpha=\alpha(\Lambda)>0$, the system
\eqref{eq} is \emph{uniformly stable} if there exists $C>0$ such
that for any $u\in H^1(\mathbb{R}_+)$ and $\Lambda\in\mathcal{P}$
we have the a priori estimate
\begin{equation}\label{stabeq}
\begin{aligned}
\alpha \|u\|^2+ |u(0)|^2 \le C(\|L(\Lambda)u\|^2/\alpha +
|\Gamma(\Lambda)u(0)|^2),
\end{aligned}
\end{equation}
where $\|\cdot\|$ denotes the $L^2(x)$ norm and $|\cdot|$ the norm
in $\mathbb{C}^k$.
\end{defi}

\begin{defi}\label{lop}
System \eqref{eq} satisfies the {\it uniform Lopatinski condition}
if
\begin{align}\label{lopeq}
\begin{split}
&(i) k=\rank\;\Gamma(\Lambda) = \dim \EE_-(\Lambda)\, \\
&(ii) |v|\le C|\Gamma(\Lambda) v| \,\text{\rm for } v\in
\EE_-(\Lambda),
\end{split}
\end{align}
for some $C>0$ independent of $\Lambda\in\mathcal{P}$.
\end{defi}

\begin{lem}\label{lop2}
Condition (ii) of the uniform Lopatinski condition \eqref{lopeq}
is a necessary and sufficient condition for $L^2(x)$ solutions of
$L(\Lambda)u=0$ to satisfy the trace estimate
\begin{equation}\label{2eq}
|u(0)|^2 \le C|\Gamma(\Lambda)u(0)|^2
\end{equation}
with $C$ independent of $\Lambda\in\mathcal{P}$.
\end{lem}

\begin{proof}
The $L^2$ solutions $u(x)$ of the constant coefficient problem
$L(\Lambda)u=0$ are precisely the functions
\begin{align}\label{a10}
u(x)=e^{xG(\Lambda)}u_0,
\end{align}
where $u_0\in E_-(\Lambda)$.

\end{proof}

The key assumption on $G(\Lambda)$ is the following one:

\begin{ass}\label{some}
For \emph{some} constant $k\times n$ matrix $\tilde \Gamma$
the system
\begin{equation}\label{tildeeq}
\begin{aligned}
u'-G(\Lambda)u&=f,\\
\tilde \Gamma u(0)&=g.
\end{aligned}
\end{equation}
has a unique solution $u\in L^2(x)$ for any given $f\in L^2(x)$,
$g\in\mathbb{C}^k$, and $u$ satisfies
\begin{align}
\alpha \|u\|^2+ |u(0)|^2 \le C(\|f\|^2/\alpha + |g|^2)
\end{align}
with $C$ independent of $f$, $g$, and $\Lambda\in\mathcal{P}$.

\end{ass}

\begin{exam}\label{symm}
It follows from the discussion in the introduction that this
assumption is satisfied by any $G(\Lambda)$ obtained as in
\eqref{ieq} by Laplace-Fourier transform of a Friedrichs
symmetrizable system \eqref{hyp}(a).    In this case $\Lambda$ and
$\mathcal{P}$ are defined as in \eqref{aa3},
$\alpha(\Lambda)=\gamma$, and we take $\tilde{\Gamma}$ to be  a
maximally dissipative boundary condition as described in Remark
\ref{aa0}, part 3.

\end{exam}

We will prove Theorem \ref{hypmain} using the following general
principle together with Example \ref{symm}.  The idea is that
existence of a boundary condition for which good estimates hold
already encodes structural properties relevant to the stability
analysis.

\begin {prop}\label{main}
Consider the resolvent-type problem $(L(\Lambda),\Gamma(\Lambda))$
as in \eqref{eq}, and suppose that $G(\Lambda)$ satisfies
Assumption \ref{some} for some choice of $\alpha(\Lambda)$.  If
the system $(L(\Lambda),\Gamma(\Lambda))$ satisfies the uniform
Lopatinski condition, then it is uniformly stable relative to
$\alpha(\Lambda)$.

\end{prop}

Proposition \ref{main} and its corollary Theorem \ref{hypmain},
proved below, extend and greatly simplify the results of
\cite{MZ2} for constant-coefficient symmetrizable systems.

\begin{proof}[Proof of Proposition \ref{main}]

Let $u\in H^1(\mathbb{R}_+)$, set
\begin{align}
L(\Lambda)u:=f, \;\;\;\Gamma(\Lambda)u(0):=g
\end{align}
and for $\tilde{\Gamma}$ as in Assumption \ref{some} introduce the
auxiliary problem
\begin{equation}\label{weq}
\begin{aligned}
w'-G(\Lambda)w&=f,\\
\tilde \Gamma w(0)&=0.
\end{aligned}
\end{equation}
By Assumption \eqref{some}, there exists a unique $L^2$ solution
$w$ satisfying
\begin{equation}\label{3}
\alpha \|w\|^2 +|w(0)|^2\le \tilde C\|f\|^2/\alpha.
\end{equation}

Now consider the residual $e:=u-w\in L^2$, satisfying
\begin{equation}\label{reseq}
\begin{aligned}
e'-G(\Lambda)e&=0,\\
\Gamma(\Lambda) e(0)&=\Gamma(\Lambda)(u(0)-w(0))= g-
\Gamma(\Lambda)(w(0)).
\end{aligned}
\end{equation}
By the uniform Lopatinski assumption and Lemma \ref{lop2},
\begin{equation}\label{eq7}
\begin{aligned}
|e(0)|^2 &\le C|\Gamma(\Lambda) e(0)|^2\\
&\le C(|g|+ |\Gamma(\Lambda) w(0)|)^2\\
& \le 2C(|g|^2+ C_1 \tilde C\|f\|^2/\alpha ).\\
\end{aligned}
\end{equation}

On the other hand, we may equally well consider \eqref{reseq} as
\begin{equation}\label{tildereseq}
\begin{aligned}
e'-G(\Lambda)e&=0,\\
\tilde \Gamma e(0)&=: \tilde g.
\end{aligned}
\end{equation}
Applying Assumption \ref{some} again, we thus have
\begin{equation}\label{eq8}
\begin{aligned}
\alpha \|e\|^2 + |e(0)|^2 &\le \tilde C|\tilde g|^2\\
&= \tilde C |\tilde \Gamma e(0)|^2,\\
\end{aligned}
\end{equation}
which, by \eqref{eq7}, gives
\begin{equation}\label{eq9}
\begin{aligned}
\alpha \|e\|^2 + |e(0)|^2 &\le
2\tilde C C_2 C(|g|^2+ C_1 \tilde C\|f\|^2/\alpha ),\\
\end{aligned}
\end{equation}
where $C_1$ is the matrix norm of $\Gamma$ and $C_2$ of $\tilde
\Gamma$. Adding \eqref{3} and \eqref{eq9}, we obtain the result
\end{proof}

\begin{proof}[Proof of Theorem \ref{hypmain}]
It remains to prove the sufficiency of the uniform Lopatinski
condition.

Let $u(x,t)\in e^{\gamma t}H^1(\mathbb{R}^{d+1}_+)$ and set
\begin{align}
Lu:=f\;\;\;\;\Gamma_\gamma u(0):=g.
\end{align}
The strong $L^2$ well-posedness of the system $(L,\Gamma_\gamma)$
follows by standard arguments (e.g., \cite{CP}, Chapter 7) from an
a priori estimate of the form \eqref{uni}  for $(L,\Gamma_\gamma)$
and an analogous estimate for the adjoint problem
$(L^*,\Gamma^*_\gamma)$.  For the definition of the adjoint
boundary condition and the verification that the adjoint problem
necessarily satisfies the (backward) uniform Lopatinski condition
provided the forward problem satisfies the (forward) Lopatinski
condition (Definition \ref{a3}), we refer to \cite{CP}, Chapter 7.

The forward estimate is an immediate consequence of  Proposition
\ref{main}, Example \ref{symm}, and Plancherel's Theorem.  The
backward estimate follows by a parallel argument, since as we
noted in Remark \ref{aa0}, if $(L,\tilde{\Gamma})$ is a
symmetrizable, maximally dissipative problem, then
$(L^*,\tilde{\Gamma}^*)$ is symmetrizable and maximally
dissipative in the backward sense. Thus, the $G(\Lambda)$ matrix
that appears in the Laplace-Fourier transformed adjoint problem
satisfies Assumption \ref{some}, and Proposition \ref{main} can be
applied to that problem as well.
\end{proof}

\bigbreak
\section{The variable-coefficient case}\label{var}

For the study of nonlinear hyperbolic boundary-value problems, it
is important to treat also the variable-coefficient analog of
\eqref{hyp},
\begin{equation}\label{varhyp}
\begin{aligned}
&L(t,x,\partial_t,\partial_x)u:=u_t + \sum_{j=1}^d A^j(t,x) u_{x_j} = f,\\
&\Gamma(t,\tilde x)
u(t,\tilde x,0)=g,\\
\end{aligned}
\end{equation}
where $L$ is Friedrichs symmetrizable and $\Gamma(t,\tilde{x})$ is
a $k\times n$ matrix or, more generally, a pseudodifferential
operator $\Gamma_\gamma(t,\tilde{x},D_t,D_{\tilde{x}})$ of degree
zero.


Strong $L^2$ well-posedness is defined for \eqref{varhyp} as in
the constant coefficient case. Following Kreiss \cite{K} we define
the uniform Lopatinski condition for \eqref{varhyp} as uniform
Lopatinski for the family of frozen-coefficient problems
\begin{equation}\label{frozenhyp}
\begin{aligned}
&u_t + \sum_{j=1}^d A^j(q,0) u_{x_j} = f\\
&\Gamma(q) u(t,\tilde x,0)=g,
\end{aligned}
\end{equation}
with parameter $q=(t,\tilde x)$ varying in $\mathbb{R}^d$, where
the constant $C>0$ is now required to be uniform in both
$\gamma>0$ and the parameter $q$. The variable-coefficient
analogue of Theorem \ref{hypmain}, extending the result proved in
\cite{K,CP,MZ2} for hyperbolic constant-multiplicity systems,
would be as follows:

The system $(L,\Gamma)$ \eqref{varhyp} is strongly $L^2$
well-posed if and only if it satisfies the frozen uniform
Lopatinski condition.

So far we have been unable to prove the sufficiency of the uniform
Lopatinski condition.  The main obstacle, curiously, is to obtain
a variable-coefficient analogue of the elementary Lemma
\ref{lop2}.  More precisely, we would like to show that if
$(L,\Gamma)$ satisfies the frozen uniform Lopatinski condition,
then solutions $u\in L^2(\mathbb{R}^{d+1}_+)$ of
\begin{align}
L(t,x,\partial_t,\partial_x)u=0
\end{align}
satisfy uniform trace estimates
\begin{align}\label{a20}
|u(t,\tilde{x},0)|_{L^2(\mathbb{R}^d)}^2\le C|\Gamma
u(t,\tilde{x},0)|_{L^2(\mathbb{R}^d)}^2.
\end{align}
For $L$ as in \eqref{varhyp} one can always find a boundary
condition $\tilde\Gamma(t,x)$ for which $(L,\tilde\Gamma)$ is
strongly $L^2$ well-posed (as in Remark \ref{aa0}, part 3),  so if
we had \eqref{a20} we could work in the original $(t,x)$ variables
and simply repeat the argument of Proposition \ref{main}, with
$\|\cdot\|_{L^2(\mathbb{R}^{d+1}_+)}$ replacing
$\|\cdot\|_{L^2(\mathbb{R}_+)}$ now, to derive the needed a priori
estimates.   In fact, in place of \eqref{a20} it would be
sufficient to establish
$$
|u(t,\tilde{x},0)|_{L^2(\mathbb{R}^d)}^2\le C\left(|\Gamma
u(t,\tilde{x},0)|_{L^2(\mathbb{R}^d)}^2+\|u\|_{L^2(\mathbb{R}^{d+1}_+)}\right).
$$
However, to do this using the tools available appears to be as
difficult as finding an actual Kreiss symmetrizer, yielding the
full estimate for general $f$. That is, the exact computation of
Lemma \ref{lop2} does not seem to be robust under lower-order
perturbations: there is no apparent advantage to small $f$ over
the general case.
%

\bigbreak
\section{Viscous shock and boundary layers}\label{viscous}
\textbf{}

In this final section we sketch how the general principle of
Proposition \ref{main} can be applied in a parabolic (or partially
parabolic) problem.

In the study of noncharacteristic viscous shock or boundary
layers,  one linearizes the compressible Navier-Stokes equations
about a function of one variable, say $w(x_d)$, which describes
the shock or boundary layer.  After symmetrizing and applying a
conjugating transformation to remove dependence on the variable
$x_d$ in the coefficients (see, e.g., the introduction to
\cite{GMWZ1} or \cite{GMWZ4}), we reduce to the study of a
constant coefficient, second-order, boundary value problem on the
half-space $\mathbb{R}^{d+1}_+$,
\begin{equation}\label{visc}
\begin{aligned}
A^0u_t + \sum_{j=1}^d A^j u_{x_j}- \sum_{j,k=1}^d B^{jk} u_{x_j,
x_k}
 &= f,\\
\Gamma u(t,\tilde{x},0)=g,
\end{aligned}
\end{equation}
where
$$
u=\begin{pmatrix} u_1\\u_2
\end{pmatrix},
\quad A^j=\begin{pmatrix}
A^j_{11} & A^j_{12} \\
A^j_{21} & A^j_{22} \\
\end{pmatrix},
\quad B^{jk}=\begin{pmatrix}
0 & 0 \\
0  & B^{jk}_{22} \\
\end{pmatrix},
$$
with $\det A^d\ne 0$, $A^0$ positive definite, $A^j$ symmetric,
and
$$
\Re \sum_{jk}\xi_j\xi_k B^{jk}_{22} \ge \theta |\xi|^2.
$$

Applying as before the Laplace transform in the temporal variable
$t$ and the Fourier transform in $\tilde x$,  we obtain the
generalized resolvent equation (with hats dropped)
\begin{equation}\label{2res}
\begin{aligned}
\lambda A^0u + &A^d u'+ \sum_{j=1}^{d-1} i\eta_j A^j u\\
&\quad -B^{dd}u'' -\sum_{j=1}^{d-1} i\eta_j (B^{j1}+B^{1j}) u'
+\sum_{j,k=1}^{d-1} \eta_j\eta_k B^{jk} u
 &= f,\\
\Gamma u(0)&=g,\\
\end{aligned}
\end{equation}
which may be written after some rearrangement as a first-order
system with a redefined $\Gamma$
\begin{equation}\label{1res}
\begin{aligned}
U'-{\CalG}(\Lambda) U&= F,\\
\Gamma U(0)&=  G,
\end{aligned}
\end{equation}
in the variable $U:=(u, u_2')$.
%

Taking ``dissipative'' boundary conditions in the class identified
by Rousset \cite{R3}, $\tilde{\Gamma} U=(\Gamma_1 u_1, u_2)$, with
$\Gamma_1$ maximally dissipative for the hyperbolic problem
$A^0_{11}v_t + \sum_j A^j_{11} v_{x_j}=0$, we obtain by
integration by parts (after forming the $L^2$ inner product of $u$
with \eqref{2res}) estimates that are nearly of the form
\eqref{stabeq}.  The difference is that several weights
$\alpha_k(\gamma, \tau, \eta)$ appear and $u$ and $u_2'$
coordinates are weighted differently.    We use this estimate to
define \emph{uniform viscous stability}, the analogue of
Definition \ref{stab}.  For this choice of weights and
$\tilde{\Gamma}$, Assumption \ref{some} is then satisfied for
$\mathcal{G}$ as in \eqref{1res}.   The explicit estimates/weights
are given in  \cite{GMWZ4,GMWZ5}.

%


A review of the proof of Proposition \ref{main} reveals that the
new weights do not affect the arguments there. Thus,  we obtain
the analogous result that uniform viscous stability follows from
the uniform Lopatinski condition.    The latter condition is
called in the viscous context the {\it uniform Evans condition}.
This extends results of \cite{MZ2, GMWZ6} in the
variable-multiplicity case, in particular for MHD.

Unfortunately, this result, though suggestive, does not yield
nonlinear stability, either for small viscosity, which requires
variable-coefficient estimates, or for large time, which requires
$L^1\to L^2$ estimates between norms \cite{GMWZ1}.

%
\begin{rem}\label{hfrmk}
\textup{ A finer point of the analysis is that the conjugating
transformation yields uniform estimates only for a \emph{compact}
set of frequencies, so a different analysis must be used in the
high-frequency regime, as discussed in \cite{GMWZ4,GMWZ5,GMWZ6}.
In particular, the Evans condition must be required to hold
uniformly under an appropriate high-frequency rescaling. However,
this high-frequency part of the analysis has already been carried
out in \cite{GMWZ4,GMWZ6} without any assumptions on multiplicity
of hyperbolic characteristics.  Thus, the bounded-frequency
argument just presented is precisely what is needed to extend to
the general, variable-multiplicity case. }
\end{rem}

\bigbreak

\end{document}